\documentclass[10pt]{article}
\usepackage{amsmath,color,amssymb,graphicx,wasysym,amsthm,setspace,ulem
}
\definecolor{darkblue}{rgb}{0.0,0.0,0.3}
\usepackage[sort&compress,numbers]{natbib}

\begin{document}
\newcommand{\eqdef}{{\ \stackrel{\mathrm{def}}{=}\ }}

\newcommand{\ra}{\rightarrow}
\newcommand{\RR}{\mathbb{R}}
\newcommand{\Sn}{\mathbb{S}^{n-1}}
\newcommand{\NN}{\mathbb{N}}
\newcommand{\RRn}{\RR^n}
\newcommand{\A}{\mathcal{A}}
\newcommand{\RRd}{\RR^n}
\newcommand{\ud}{\mathrm{d}}
\newtheorem{theorem}{Theorem}[section]
\newtheorem{lemma}{Lemma}[section]
\newtheorem{example}[theorem]{Example}
\newtheorem{remark}[theorem]{Remark}
\bibliographystyle{plain}
\normalem

\title{Structured and unstructured continuous models for \textit{Wolbachia} infections} 
\author{J\'ozsef Z. Farkas\thanks{Department of Computing Science and
Mathematics, University of Stirling, FK9 4LA, Scotland, UK; email:
\texttt{jzf@maths.stir.ac.uk}} \and Peter Hinow\thanks{Department of
Mathematical Sciences, University of Wisconsin -- Milwaukee, P.O.~Box 413,
Milwaukee, WI 53201, USA; phone: +1 414 229 4933; email:
\texttt{hinow@uwm.edu}}} 
\date{\today}
\maketitle

\begin{abstract}
We introduce and investigate a series of models for an infection of a diplodiploid host species by the bacterial endosymbiont \textit{Wolbachia}. The continuous models are characterized by partial vertical transmission, cytoplasmic incompatibility and fitness costs associated with the infection. A particular aspect of interest is competitions between mutually incompatible strains. We further introduce an age-structured model that takes into account different fertility and mortality rates at different stages of the life cycle of the individuals. With only a few parameters, the ordinary differential equation models exhibit already interesting dynamics and can be used to predict criteria under which a strain of bacteria is able to invade a population. Interestingly, but not surprisingly, the age-structured model shows significant differences concerning the existence and stability of equilibrium solutions compared to the unstructured model.

\begin{sloppypar}
\end{sloppypar}
\end{abstract}

\section{Introduction}\label{section:Introduction}

\textit{Wolbachia} is a maternally transmitted bacterium that lives in symbiosis
with many arthropod species and some filarial nematodes
\citep{Werren1997,ONeill}. It inhabits testes and ovaries of its hosts and has
the ability to interfere with their reproductive mechanisms, resulting in a
variety of phenotypes. Well known effects are cytoplasmic incompatibility,
induction of parthenogenesis, and feminization of genetic males, depending on
the host species and the \textit{Wolbachia} type. Besides the intrinsic interest
in these mechanisms, \textit{Wolbachia} are investigated as tools to drive
desirable genes into a target population \citep{Rasgon2004}, as reinforcers of
speciation  \citep{Telschow2005,Telschow2005b,Keeling2003}, and as potential
means of biological control  \citep{McMeniman2009}. It was recently shown by
McMeniman \textit{et al.}~\citep{McMeniman2009} that infection with
\textit{Wolbachia} shortens the life span of the mosquito \textit{Aedes
aegypti}, a vector for the Dengue fever virus. Since only older mosquitoes are
carriers, this is a promising strategy to reduce the transmission of pathogens,
without the ethically untenable eradication of a vector species.

Beginning already a half a century ago \citep{Caspari1959}, various mathematical
models for the spread of a  \textit{Wolbachia} infection have been proposed and
studied in the literature, see
e.g.~\citep{Turelli1994,Rasgon2004,Telschow2005,Keeling2003,Engelstadter2004,
Schofield2002,Vautrin2007,Haygood2009} and references therein. Largely, these
models fall into two classes, depending on whether time proceeds in discrete
steps or continuously. Examples for continuous  models are the papers
\citep{Keeling2003} and \citep{Schofield2002} that employ ordinary, respectively
partial differential equations (with a spatial structure in the latter case). In
their paper \citep{Keeling2003}, the authors proposed and studied a simple
continuous model for the infection of an arthropod population with cytoplasmic
incompatibility (CI) causing \textit{Wolbachia}. Cytoplasmic incompatibility in
diplodiploid (i.e.~with diploid males and females) species manifests itself in
completely or partially unviable crosses of infected males with uninfected
females. For a discussion of the more complex outcome of cytoplasmic
incompatibility in haplodiploid species see \citep{Vautrin2007}.

In this paper, we introduce a series of models for different aspects of interest. We start in section \ref{section:single_sex_model} with an ordinary differential equation model for a single \textit{Wolbachia} strain that infects a population without separate sexes. In section \ref{section:multiple_strains} we present a model for infections with multiple strains. The present theoretical literature offers a complex picture of infection with multiple strains. While some authors exclude the coexistence of multiple strains of \textit{Wolbachia} in infected individuals \cite{Keeling2003,Haygood2009}, others model doubly infected individuals as a class of their own \cite{Engelstadter2004,Vautrin2007}. Moreover, different assumptions are made about the mutual compatibility of individuals carrying different strains. We construct a general model that encompasses these different possibilities by suitable choices of parameter values and/or initial conditions. Finally, motivated by the study \cite{McMeniman2009}, in section \ref{section:age_structure} we  refine our model from section \ref{section:single_sex_model}, by considering age-structured populations.
We refer the interested reader for basic concepts and results in structured
population dynamics to \citep{Cushing, Metz, Webb}. Modeling structured
populations usually involves partial differential equations which are more    
difficult to analyze. Analytical progress is still possible, and as we will see
in \ref{section:age_structure}, the age-structured model may exhibit richer
dynamics. At this point it will be possible to study age-dependent killing of
\textit{Wolbachia} infected individuals. Our models contain parameters of only
three types, namely transmission efficacies, levels of cytoplasmic
incompatibility and fitness costs for the infected individuals. The analysis of
the models aims to give conditions for the stability of specific equilibrium
solutions that correspond to successful invasions. The paper ends with a
discussion in section \ref{section:discussion} and an outlook on future work.

\begin{sloppypar}
\section{Single-sex model for a singular \textit{Wolbachia} strain}\label{section:single_sex_model}
\end{sloppypar}
Assume that the ratio of infected males to infected females is the same as the
ratio of uninfected males  to uninfected females and hence the population can be
formally considered hermaphroditic. Let $I$ and $U$ denote the number of
infected, respectively uninfected individuals in the population. Vertical
transmission is partial, let $\tau\in[0,1]$ be the fraction of infected
offspring from infected parents (another common notation is $\mu=1-\tau$ for the
fraction of uninfected ova produced by an infected female, see
e.g.~\citep{Turelli1994,ONeill,Vautrin2007}). Furthermore, we follow Keeling
\textit{et al.} \citep{Keeling2003} and assume that the birth rate for both
infected and uninfected individuals is equal (no reduction in fecundity in
infected individuals). Let this rate be denoted by $b>0$. Death of the
individuals is modelled by a logistic loss term with rate $d>0$ that accounts
for competition among the total population. However, infected individuals can
suffer an additional loss of fitness given by $D\ge0$. Cytoplasmic
incompatibility arises when an infected male fertilizes an egg from an
uninfected female. Then, with a probability $q\in[0,1]$, the offspring is
nonviable. As we do not consider separate sexes in this simple model, we just
reduce the amount of offspring from uninfected individuals based on the
probability of an encounter with an infected individual. Uninfected individuals
still arise due to incomplete transmission of the bacterium from infected
parents. Our equations read
\begin{align*}
\frac{dI}{dt} &= \left(\tau b -(d+D)(I+U)\right)I, \\
\frac{dU}{dt} &= (1-\tau)bI +\left(b\left(1-q\frac{I}{I+U}\right) -d(I+U)\right)U.
\end{align*}
Upon rescaling the time by $t\ra bt$ and setting $\eta=\frac{d+D}{d}$, we obtain the reduced system for the quantities $i=dI/b$, $u=dU/b$
\begin{align}
\frac{di}{dt} &= (\tau -\eta(i+u))i, \label{eq_I}\\
\frac{du}{dt} &= (1-\tau)i +\left( 1-q\frac{i}{i+u}  -(i+u)\right)u. \label{eq_U}
\end{align}
Observe that $\eta^{-1}=:\xi\in(0,1]$ and that $1-\xi$ can be interpreted as the
fitness cost associated with \textit{Wolbachia} infection. The point $(0,0)$ can
be added to the domain of the state space, with the understanding that it is an
equilibrium solution. It is obvious that the subspace of uninfected populations
$\{0\}\times\RR$ is forward invariant (that is, if initially there are no
infected individuals, then there will be none at later times) and if
transmission is complete, $\tau=1$, then so is the subspace  of completely
infected populations $\RR\times\{0\}$.  

\begin{figure}[th]
\begin{center}
\includegraphics[width=60mm]{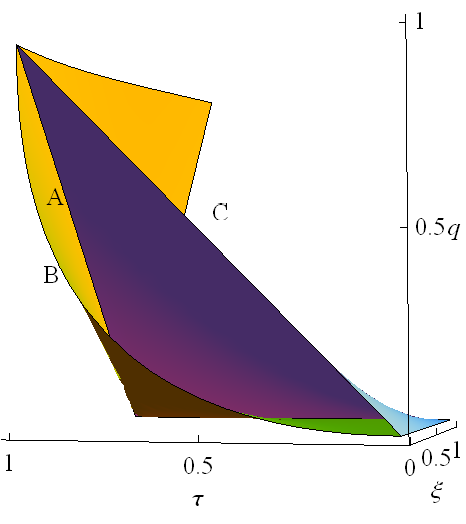}
\caption{The yellow surface separates the $(\xi,\tau,q)$-parameter space of
model \eqref{eq_I}--\eqref{eq_U} into a region where only the disease free
equilibrium \eqref{disease_free} exists (A) and where coexistence equilibrium
solutions $(u_1,i_1)$ and $(u_2,i_2)$ given by \eqref{coexist} are possible (B).
However, only in the region (C) above the blue surface  given by
\eqref{u2_non_neg} is $u_2=\tau\xi-i_2\ge0$ (this belongs to the observed
stable equilibrium $(u_2,i_2)$).
}\label{fig_1}
\end{center}
\end{figure}

Equation \eqref{eq_I}--\eqref{eq_U} always admits the disease free equilibrium
\begin{equation}\label{disease_free}
(i_0,u_0) = (0,1).
\end{equation}
Setting the left-hand side of equation \eqref{eq_I} to zero and solving for an  equilibrium point $u$ yields
\begin{equation}\label{u_from_i}
u = \tau\xi-i.
\end{equation}
Inserting this into the equilibrium condition for \eqref{eq_U} gives a quadratic equation for $i$,
\begin{equation}\label{quadratic}
\frac{q}{\tau\xi}i^2 +(\tau(\xi-1)-q )i +\tau\xi(1-\xi\tau)=0.
\end{equation}
Provided that
\begin{equation}\label{discriminant}
(\tau(\xi-1)-q )^2 -4q(1-\xi\tau)\ge0,
\end{equation}
equation \eqref{quadratic} has the solutions
\begin{equation}\label{coexist}
\begin{aligned}
i_1=i_1(q,\tau,\xi) &= \frac{\tau\xi\left( q-\tau(\xi-1) - \sqrt{(\tau(\xi-1)-q )^2 -4q(1-\xi\tau)}\right)}{2q},  \\
i_2=i_2(q,\tau,\xi) &= \frac{\tau\xi\left( q-\tau(\xi-1) + \sqrt{(\tau(\xi-1)-q )^2 -4q(1-\xi\tau)}\right)}{2q},
\end{aligned}
\end{equation}
\begin{sloppypar}\noindent
these are always non-negative. The corresponding equilibrium values for the
uninfected individuals are given by \eqref{u_from_i}. For $u_2\ge0$ it is
necessary that 
\begin{equation}\label{u2_non_neg}
q+\tau(\xi-1)\ge0. 
\end{equation}
This condition is also sufficient,
since then one can derive from $\tau\le1$ the inequality 
\end{sloppypar}
\begin{equation*}
q+\tau(\xi-1) \ge \sqrt{(\tau(\xi-1)-q )^2 -4q(1-\xi\tau)}
\end{equation*}
and hence  $u_2\ge0$.  Condition \eqref{discriminant} separates
two regions in the parameter space, depending on whether other equilibrium
solutions than the disease free equilibrium are possible (see figure
\ref{fig_1}). We calculate the Jacobian of
the right hand side $F$ of \eqref{eq_I}--\eqref{eq_U},
\begin{equation}\label{Jacobian}
DF(i,u) = \begin{pmatrix} \tau-\eta(2i+u) & -\eta i\! \\ \!1-\tau+qu\left(\frac{i}{(i+u)^2}-\frac{1}{i+u}\right)-u & 1+qi\left(\frac{u}{(i+u)^2}-\frac{1}{i+u}\right)\!-(2u+i)\end{pmatrix}\!.
\end{equation}
At the  disease free equilibrium \eqref{disease_free}, we have
\begin{equation*}
DF(0,1) = \begin{pmatrix} \tau-\eta  & 0 \\ -\tau-q  & -1 \end{pmatrix},
\end{equation*}
this matrix has the eigenvalues $-1$ and $\tau-\eta\le 0$. The latter eigenvalue
is $0$ only if $\tau=1$ (complete transmission) and $D=0$ (no penalty for
infection), in all other cases it is strictly negative, and the  disease free
equilibrium is locally asymptotically stable. 

Explicit expressions (with respect to the parameters $q,\,\tau$ and $\xi$) can
be obtained for the eigenvalues of the Jacobian \eqref{Jacobian} at the 
equilibrium solutions using e.g.~\textsc{mathematica} (the
\textsc{mathematica} notebook is available from the authors upon request).
Unfortunately these expressions are very complicated and not easily analyzed. We
will present instead some representative numerical examples to demonstrate the
possible behaviours.

\begin{figure}[th]
\begin{center}
\includegraphics[width=60mm]{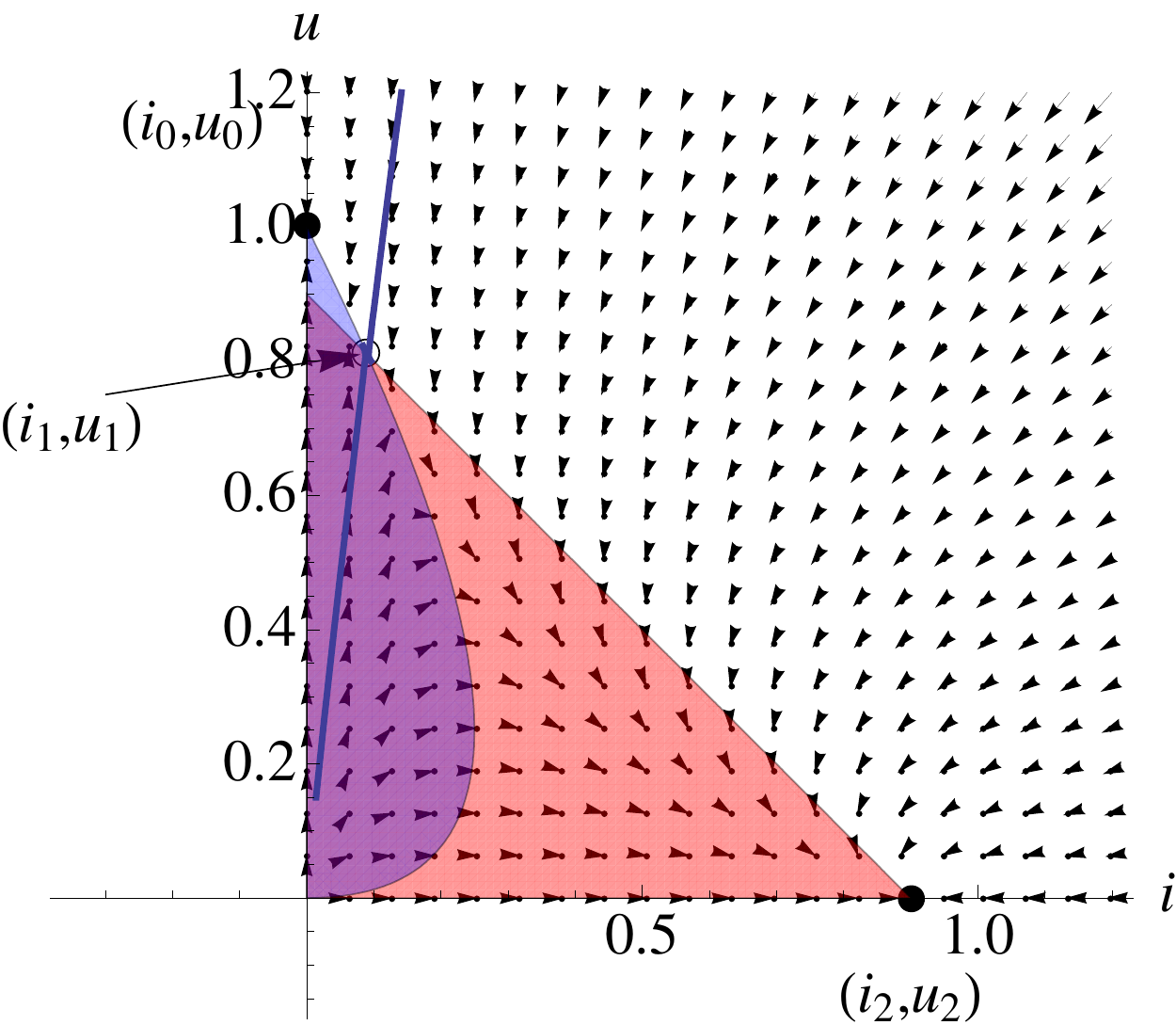}\includegraphics[width=60mm]{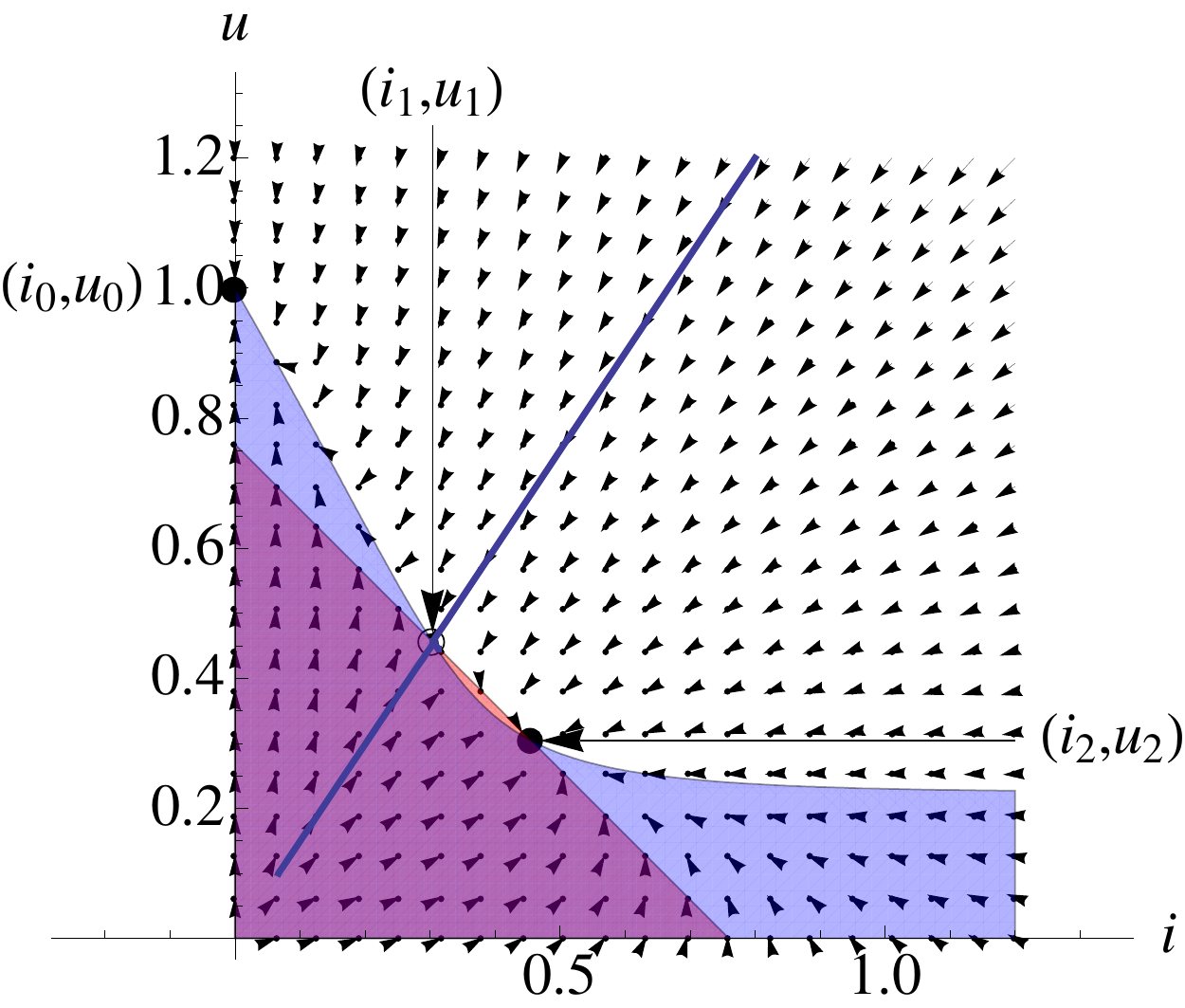}
\caption{\textit{(Left)} The vector field \eqref{eq_I}--\eqref{eq_U} for the
parameter triple $(\xi,\tau,q)=(0.9,1,1)$ together with the three equilibrium
solutions. Solid disks indicate locally asymptotically stable equilibrium
solutions, while the disk indicates an unstable equilibrium. Shown are also
regions of growth of $u$ (light blue) and growth of $i$ (light red). The blue
lines are the stable manifolds of the saddle point $(i_1,u_1)$ and the
separatrices of the equilibrium solutions $(i_0,u_0)$ and $(i_2,u_2)$.
\textit{(Right)} The vector field \eqref{eq_I}--\eqref{eq_U} for the parameter
triple $(\xi,\tau,q)=(1,0.76,1)$, which admits bistability and true coexistence.
The stable manifold of the saddle point  $(i_1,u_1)$ is shown in
blue.}\label{fig_2}
\end{center}
\end{figure}

\begin{example}\label{examp_1} {\rm Assume that the infection is completely
inherited, $\tau=1$, cytoplasmic incompatibility is complete $q=1$, and that the
cost of the infection is low, $\xi=0.9$. Then the three equilibrium solutions
$(i_0,u_0)=(0,1)$, $(i_1,u_1)=(0.09,0.81)$ and $(i_2,u_2)=(0.9,0)$ are present,
of which $(i_0,u_0)$  and $(i_2,u_2)$ are locally asymptotically stable. The
vector field is shown in figure \ref{fig_2} (left panel). The epidemic
equilibrium $(i_2,u_2)$ has a much bigger basin of attraction than the disease
free equilibrium $(i_0,u_0)$, in other words, the threshold for an infection to
establish itself is low.}
\end{example}
\begin{example} {\rm
At high levels of cytoplasmic incompatibility, $q=1$, and no penalty for the
infection, $\xi=1$, and a partial transmission $\tau=0.76$, besides the
equilibrium $(i_0,u_0)$ there is the locally asymptotically stable coexistence
equilibrium $(i_2,u_2)=(0.456,0.304)$,  see figure
\ref{fig_2} (right panel). The basin of attraction of the disease free
equilibrium $(i_0,u_0)$ is much bigger than in example \ref{examp_1}, indicating
that the threshold for establishment of the infection is higher.}
\end{example}
\begin{example} {\rm
If the cytoplasmic incompatibility is very weak, $q=0.1$ and the fitness cost of
the infection is higher, $\xi=0.5$ (and again $\tau=1$), then the equilibrium
$(i_0,u_0)=(0,1)$ is globally asymptotically stable in the set
\mbox{$\{i\ge0,u\ge0\}$} and the equilibrium $(i_1,u_1)=(0.5,0)$ is a saddle
point. The third equilibrium point has $u_2<0$ and has no biological meaning. }
\end{example}

The above examples suggest that the equilibrium $(i_1,u_1)$ is always a saddle
and that $(i_2,u_2)$ is locally asymptotically stable, if $u_2\ge 0$. If
$(i_2,u_2)$ exists and is locally asymptotically stable, then its basin of
attraction is larger for larger values of the transmission rate $\tau$ (this is
the intersection of regions (B) and (C) in figure \ref{fig_1}).

\section{Infection with multiple strains}\label{section:multiple_strains}
\begin{figure}[th]
\begin{center}
\includegraphics[width=60mm]{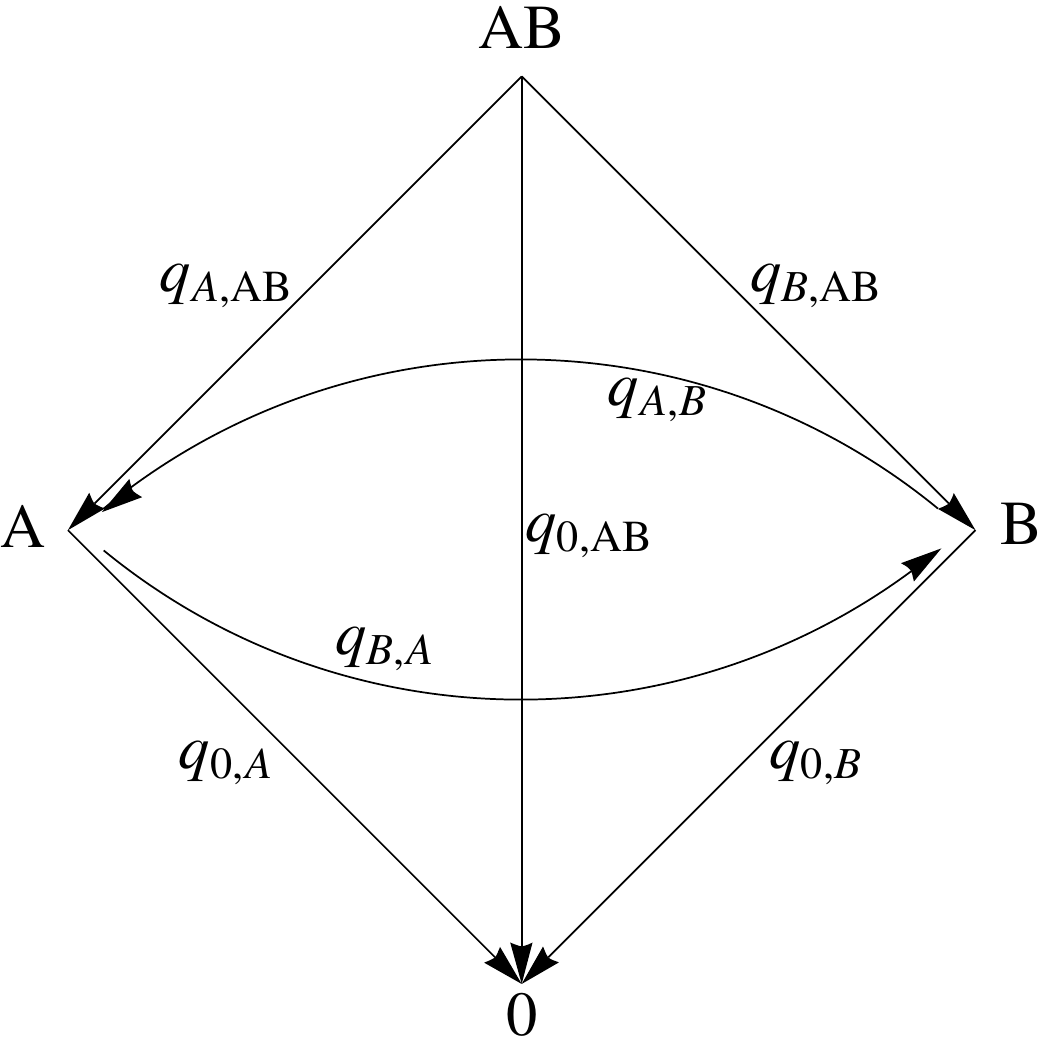}
\caption{The directed graph of possible incompatibility relations. An arrow from
node $X$ to node $Y$ indicates that a $X\male\times Y\female$ cross is
incompatible, with incompatibility level $q_{Y,X}$.}\label{fig_3}
\end{center}
\end{figure}

Assume that two \textit{Wolbachia} strains \textit{A} and \textit{B} are present
in the population and let $i_A,\,i_B$ and $i_{AB}$ denote the number of
individuals singly or doubly infected, in addition to $u$, the number of
uninfected individuals. We use the same scaling as in the previous section,
where time was scaled with respect to the birth rate $b$. There are seven
incompatibility levels $q_{X,Y}$ where $X,Y\in\{0,A,B,AB\}$ are the infection
statuses in the incompatible cross. These are shown in the directed graph of
incompatibilities in figure \ref{fig_3}. 

We assume that the transmission of one strain in doubly infected individuals is
independent of the transmission of the other strain  \citep{Vautrin2007}.
Moreover, the mortalities due to infection with either strain are additive. With
these assumptions, and setting $p=i_{AB}+i_A+i_B+u$ for the total population,
our model is
\begin{equation}\label{multiple_strains}
\begin{aligned}
\frac{di_{AB}}{dt} &= \tau_A\tau_B i_{AB} -\eta_{AB} pi_{AB}, \\
\frac{di_A}{dt}    &= \tau_A(1-\tau_B) i_{AB} + \tau_A\left(1-q_{A,B}\frac{i_B}{p}-q_{A,AB}\frac{i_{AB}}{p}\right)i_A-\eta_Api_A, \\
\frac{di_B}{dt}    &= (1-\tau_A)\tau_B i_{AB} + \tau_B\left(1-q_{B,A}\frac{i_A}{p}-q_{B,AB}\frac{i_{AB}}{p}\right)i_B-\eta_Bpi_B, \\
\frac{du}{dt}      &= (1-\tau_A)(1-\tau_B)i_{AB} + (1-\tau_A)\left(1-q_{A,B}\frac{i_B}{p}-q_{A,AB}\frac{i_{AB}}{p}\right)i_A \\
&\quad + (1-\tau_B)\left(1-q_{B,A}\frac{i_A}{p}-q_{B,AB}\frac{i_{AB}}{p}\right)i_B \\
&\quad + \left(1-q_{0,A}\frac{i_A}{p}-q_{0,B}\frac{i_B}{p}-q_{0,AB}\frac{i_{AB}}{p}\right)u -pu,
\end{aligned}
\end{equation}
where
\begin{equation*}
\eta_A = \frac{d+D_A}{d},  \quad \eta_B = \frac{d+D_B}{d} \quad \textrm{and} \quad \eta_{AB}=\frac{d+D_A+D_B}{d} =\eta_A+\eta_B-1
\end{equation*}
are measures of the fitness costs of the individual infection types. This model can now be reduced in its complexity in a variety of ways. For example, if the initial condition lies in the space $\{0\}\times\RR^3$, that is, there are no doubly infected individuals present initially, then the solution will also lie in that space at all times. Moreover, setting appropriate incompatibility levels $q_{X,Y}$ to zero allows to study cases of mutual compatibility.

As a first illustration, we want to consider the absence of doubly infected individuals, $i_{AB}\equiv 0$, mutual compatibility of infected individuals, $q_{A,B}=q_{B,A}=0$, equal transmission efficacy  $\tau_A=\tau_B=:\tau$ and equal infection costs $\eta_A=\eta_B=:\eta$ (again, we write $\xi=\eta^{-1}$). The symmetry is only broken by the levels of incompatibility $q_{0,A}\neq q_{0,B}$. Hence, equation \eqref{multiple_strains} simplifies to
\begin{equation}\label{multiple_strains_simple}
\begin{aligned}
\frac{di_A}{dt} &= (\tau -\eta p)i_A, \\
\frac{di_B}{dt} &= (\tau -\eta p)i_B, \\
\frac{du}{dt}   &= (1-\tau)(i_A + i_B)  +
\left(1-q_{0,A}\frac{i_A}{p}-q_{0,B}\frac{i_B}{p}\right)u -pu.
\end{aligned}
\end{equation}
Observe that model \eqref{multiple_strains_simple} has the property that planes  orthogonal  to the $(i_A,i_B)$-plane
\begin{equation*}
R_\alpha= \left\{(i_A,i_B,u) \in \RR_{\ge0}^3\::\: i_A-\alpha i_B =  0 \right\}
\end{equation*}
are invariant under the flow generated by  \eqref{multiple_strains_simple}. This is seen from the fact that for every $\alpha\in[0,\infty]$
\begin{equation*}
\frac{d}{dt}(i_A-\alpha i_B) = (\tau -\eta p)(i_A-\alpha i_B) = 0
\end{equation*}
on $R_\alpha$. This implies that the ratio $\frac{i_A}{i_A+i_B}$ remains
constant along a trajectory. In other words, if transmission efficacies and
death rates are equal for two strains (as are birth rates throughout in our
model) then neither strain can replace the other among the infected individuals
based on stronger cytoplasmic incompatibility. This is in line with recent
theoretical predictions of Turelli \citep{Turelli1994} and Haygood and Turelli
\citep{Haygood2009}, who suggest that strains are selected for relative
fecundity rather than high levels of cytoplasmic incompatibility. It needs to be
pointed out however, that even a small difference in transmission efficacies or
death rates of the two strains helps the strain with the greater transmission
rate or the lower mortality to establish itself in the population.

The disease free equilibrium of \eqref{multiple_strains_simple} is easily found to be
\begin{equation}\label{multi_disease_free}
(i_{A,0},i_{B,0},u_0) = (0,0,1).
\end{equation}
It is clear that the subspaces $\{0\}\times\RR\times\RR$ and  $\RR\times\{0\}\times\RR$ are forward invariant under the flow generated by \eqref{multiple_strains_simple} and that the equilibrium solutions  \eqref{coexist} exist in the respective subspaces, with $q$ in \eqref{coexist} replaced by either $q_{0,A}$ or $q_{0,B}$. That is, we have equilibrium solutions
\begin{equation*}
\begin{aligned}
(i_{A,1}(q_{0,A},\tau,\xi),0,\tau\xi-i_{A,1}), &\quad (i_{A,2}(q_{0,A},\tau,\xi),0,\tau\xi-i_{A,2}), \\
(0,i_{B,1}(q_{0,B},\tau,\xi),\tau\xi-i_{B,1}), &\quad
(i_{B,2}(0,q_{0,B},\tau,\xi),\tau\xi-i_{B,2}).
\end{aligned}
\end{equation*}
Besides that,  it can be checked that there is a
continuum of equilibrium solutions 
\begin{equation}\label{cont_equisols}
\begin{pmatrix} i_A(u) \\ i_B(u) \end{pmatrix} =
\begin{pmatrix}
\frac{-\tau^3-\eta q_{0,B}\tau u+\eta^2 q_{0,B}
u^2+\tau^2(1+(\eta-1)u)}{\eta^2(q_{0,A}-q_{0,B}) u} \\
\frac{\tau^3+\eta q_{0,A}\tau u-\eta^2 q_{0,A}
u^2-\tau^2(1+(\eta-1)u)}{\eta^2(q_{0,A}-q_{0,B}) u}
\end{pmatrix}
\end{equation}
for every $u\in(0,\tau\xi)$, provided that these expressions are  non-negative. The solutions from \eqref{cont_equisols}  satisfy
\begin{equation*}
i_A(u) + i_B(u) = \frac{\tau}{\eta}-u = \tau\xi-u,
\end{equation*}
which corresponds to equation \eqref{u_from_i}.

\begin{figure}[th]
\begin{center}
\includegraphics[width=60mm]{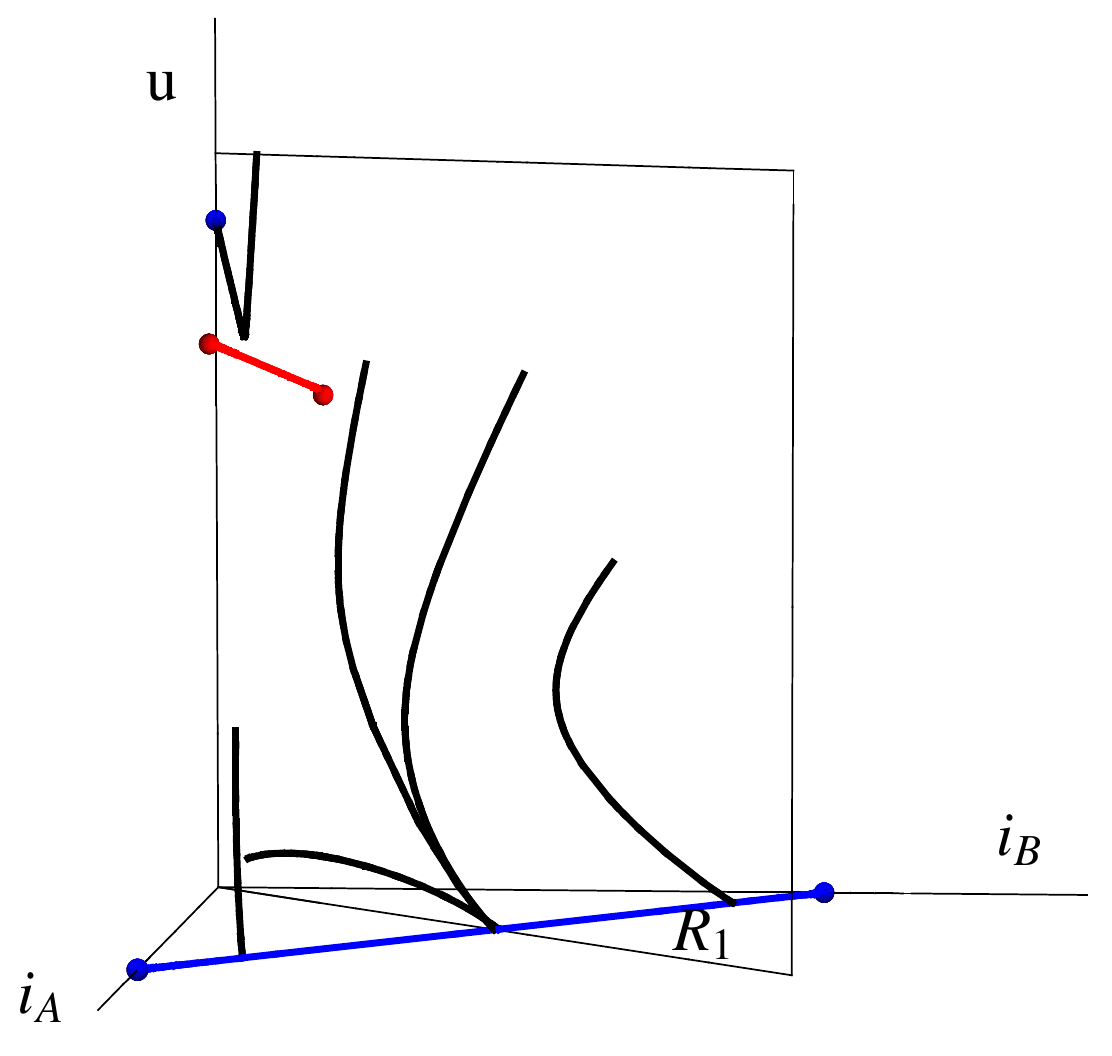}\includegraphics[width=60mm]{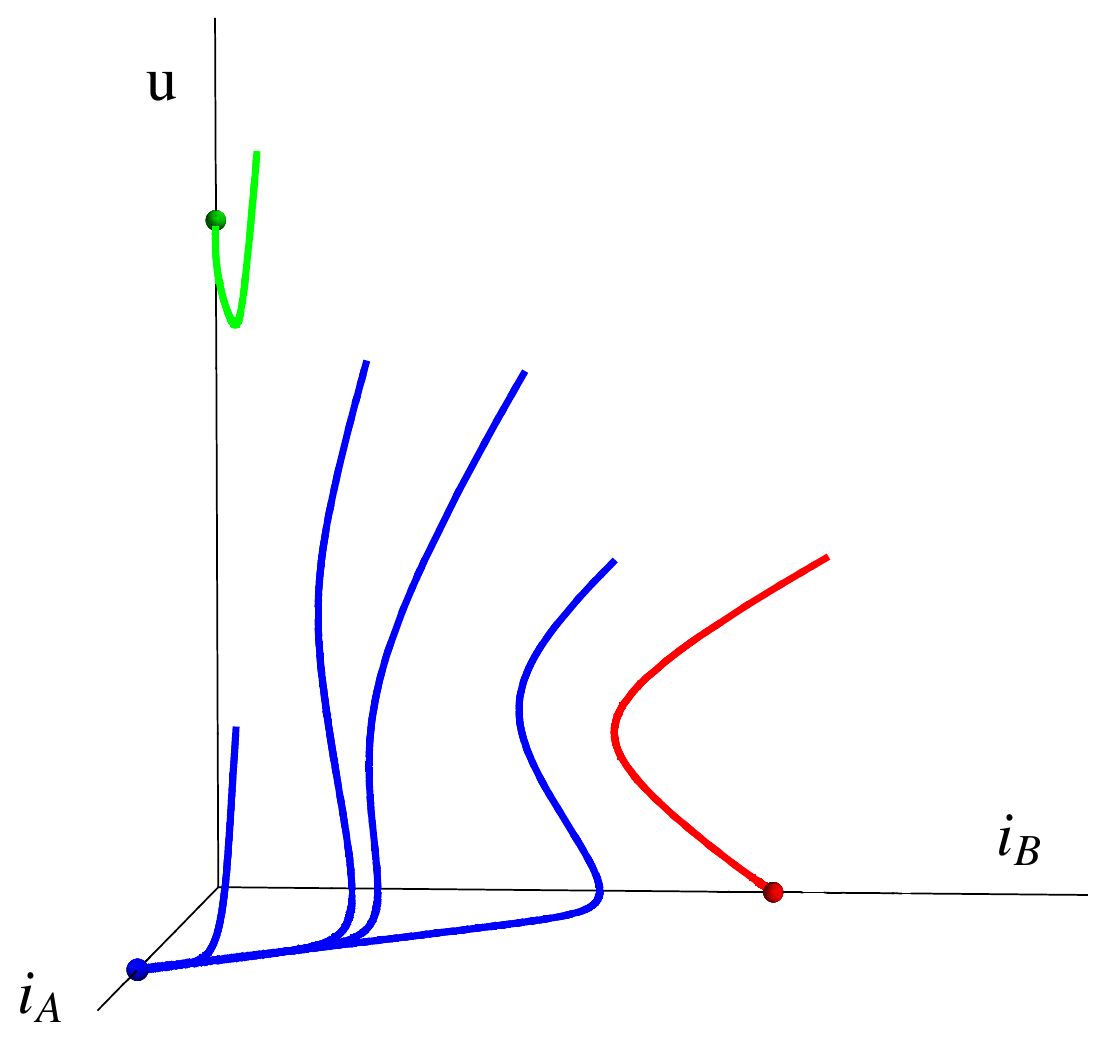}
\caption{\textit{(Left)} Dynamics of model \eqref{multiple_strains_simple} in
the case that $q_{0,A}=0.95,\,q_{0,B}=0.5,\,\tau=1$ and $\eta=1.1$. The solid
blue line is a family of attractors for the flow restricted to planes
$R_\alpha$, the red line  is a branch of saddle points in each of these
subspaces. Several individual trajectories are shown, and the space
\mbox{$R_1=\{i_A=i_B\}$} is marked. \textit{(Right)} If we choose instead
$\eta_A=1.1$ and $\eta_B=1.2$ (while keeping all other parameters the same) then
the less costly strain \textit{A} dominates the population.}\label{fig_4}
\end{center}
\end{figure}

\begin{example} {\rm We consider the case $\tau=1$ of complete transmission. One checks by direct calculation that system \eqref{multiple_strains_simple} then has another manifold of equilibrium solutions
\begin{equation*}
i_A  + i_B  = \xi, \quad u = 0.
\end{equation*}
The intersection of this line with each plane $R_\alpha$ orthogonal to the 
$(i_A,i_B)$-plane is an equilibrium for the flow restricted to that plane. 
In
addition, each $R_\alpha$ contains a saddle point. A numerical example is shown
in figure \ref{fig_4} (left panel). If, in contrast, different costs are
associated with the infection, the strain with the lower cost will dominate the
population, see figure \ref{fig_4} (right panel). Similarly, if infection costs
are equal, but one strain transmits  more efficiently, then it is going to
dominate the population.}
\end{example}

We return to model \eqref{multiple_strains} and consider mutual incompatibility
of singly infected individuals. The equations are
\begin{equation}\label{mutually_incompatible}
\begin{aligned}
\frac{di_A}{dt} &= \left(\tau\left(1-q_{A,B}\frac{i_B}{p}\right) -\eta p\right)i_A, \\
\frac{di_B}{dt} &= \left(\tau\left(1-q_{B,A}\frac{i_A}{p}\right) -\eta p\right)i_B, \\
\frac{du}{dt}   &= \phantom{+}(1-\tau)\left(\left(1-q_{A,B}\frac{i_B}{p}\right)i_A+\left(1-q_{B,A}\frac{i_A}{p}\right)i_B\right) \\
&\quad +\left(1-q_0\frac{i_A+i_B}{p}\right)u-pu.
\end{aligned}
\end{equation}
It is clear that if either strain is not present initially then it will remain
absent at all times. On the marginal spaces $\{i_A=0\}$ and $\{i_B=0\}$,
equations \eqref{mutually_incompatible} reduce to the single strain model
\eqref{eq_I}--\eqref{eq_U} and have the corresponding equilibrium solutions
where only one strain is present, with the common disease free equilibrium
\eqref{multi_disease_free}. It follows from \eqref{mutually_incompatible} that

\begin{equation*}
\begin{aligned}
\frac{d}{dt} (q_{B,A}i_A-q_{A,B}i_B) &= \phantom{-}q_{B,A}\left(\tau\left(1-q_{A,B}\frac{i_B}{p}\right)-\eta p\right)i_A \\
&\quad -q_{A,B}\left(\tau\left(1-q_{B,A}\frac{i_A}{p}\right) -\eta p\right)i_B \\
&=\phantom{+}\tau\left(q_{B,A}\left(1-q_{A,B}\frac{i_B}{p}\right)i_A-q_{A,B}
\left(1-q_{B,A}\frac{i_A}{p}\right)i_B\right)\\
&\quad  - \eta p(q_{B,A}i_A-q_{A,B}i_B) \\
&=\left(\tau-\eta p\right)(q_{B,A}i_A-q_{A,B}i_B).
\end{aligned}
\end{equation*}
This implies that the plane orthogonal to the $(i_A,i_B)$-plane,
\begin{equation*}
R_*= \left\{(i_A,i_B,u) \in \RR_{\ge0}^3\::\: q_{B,A}i_A = q_{A,B}i_B  \right\},
\end{equation*}
is forward invariant, and hence so are the wedges on either side. Solving
equations \eqref{mutually_incompatible} for the total population yields
that if $i_A\neq 0$, $i_B\neq 0$ 
\begin{equation*}
\begin{aligned}
0 &= p^2-\frac{\tau}{\eta}p + \frac{\tau}{\eta}q_{A,B}i_B, \\
0 &= p^2-\frac{\tau}{\eta}p + \frac{\tau}{\eta}q_{B,A}i_A.
\end{aligned}
\end{equation*}
For this system to be consistent, it is necessary that if $i_A\neq 0$, $i_B\neq
0$ then
\begin{equation*}
q_{B,A}i_A = q_{A,B}i_B,
\end{equation*}
hence any coexistence equilibrium of the two infected strains has to lie in the
plane $R_*$. Indeed, there may be a saddle point equilibrium solution in $R_*$
that is locally asymptotically stable for flows that start in $R_*$.

\begin{figure}[th]
\begin{center}
\includegraphics[width=60mm]{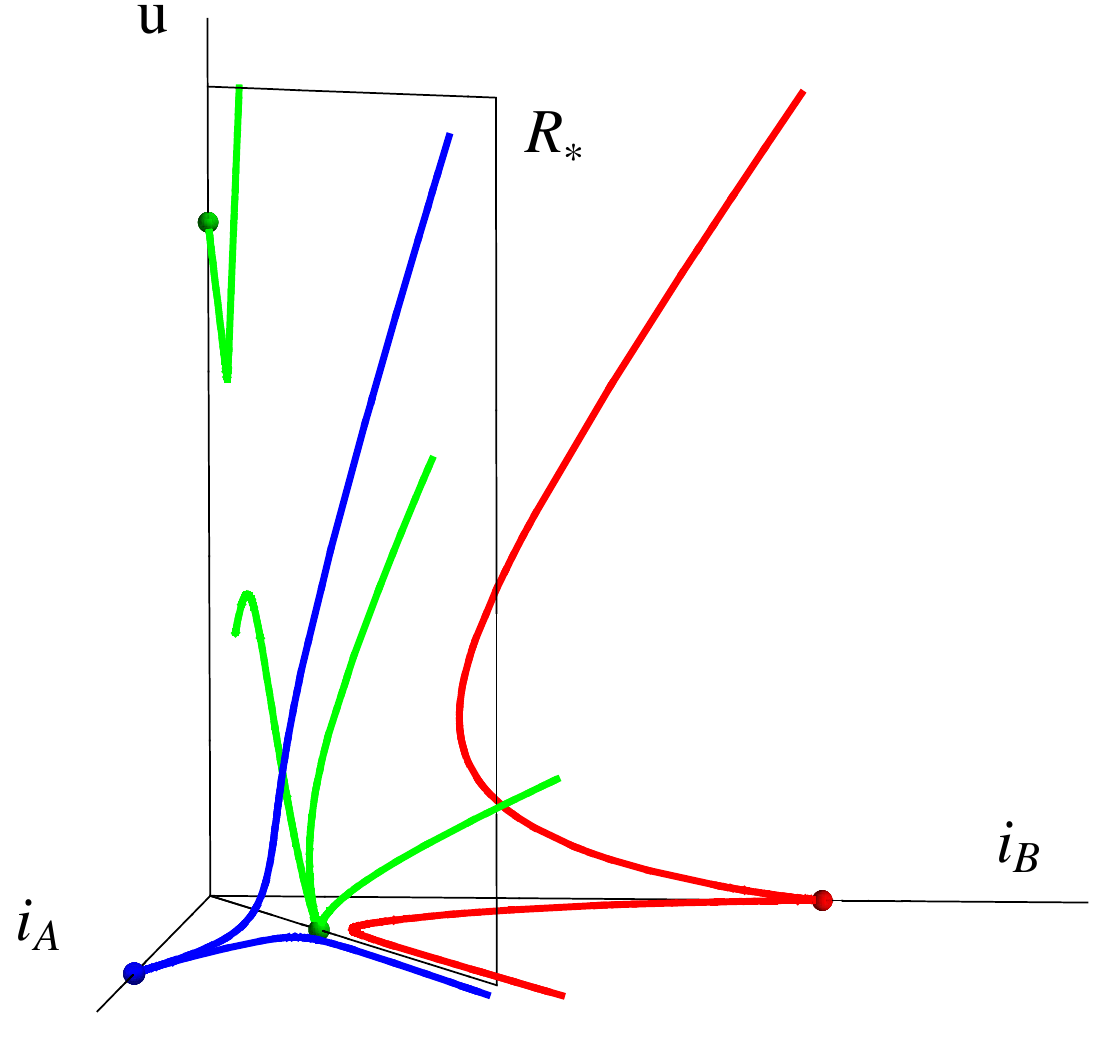}
\caption{Mutually incompatible strains as described in
\eqref{mutually_incompatible} do not show coexistence. The marginal equilibrium
solutions are marked by dots, the disease free equilibrium (equation 
\eqref{multi_disease_free}, green) attracts some solutions from the plane $R_*$
and
other trajectories beginning in the wedges on either side that
are not shown here.}\label{fig_5}
\end{center}
\end{figure}

\begin{example} {\rm Let $q_{A,B}=0.99$, $q_{B,A}$ and $q_{0,A}=q_{0,B} = 1$.
Choose the transmission efficacy $\tau=1$ and the cost of the infection
$\eta=1.1$. We see in figure \ref{fig_5} that every solution starting
outside the space $R_*$ converges to an equilibrium in one of the marginal
spaces. The plane $R_*$ contains a locally stable equilibrium for trajectories
starting in $R_*$  which has $u_2=0$ and
it contains a saddle point for trajectories
starting within $R_*$ (not shown, compare to example \ref{examp_1} and figure
\ref{fig_2}).
}
\end{example}
Again, we need to recall that this dynamical behaviour is not generic, in the
sense that the complement of the set $\{\tau_A=\tau_B,\,\eta_A=\eta_B \}$ is
dense in the parameter space.

Finally, we want to explore the full model \eqref{multiple_strains} when doubly infected individuals are present. To reduce  the number of parameters somewhat, we assume
\begin{equation*}
q_{A,B} = q_{A,AB} = q_{0,B}, \quad q_{B,A} = q_{B,AB} = q_{0,A},
\end{equation*}
that is, the presence of one strain in the fertilizing male that is missing in
the female has the same effect, regardless of the other infections that the
female may carry \citep[p.66]{Hoffmann1997}. Moreover, the escape from
cytoplasmic incompatibility for the offspring of an uninfected female and a
doubly infected male is the product of the two individual escape probabilities
\citep{Vautrin2007},
\begin{equation*}
1-q_{0,AB} = (1-q_{0,A})(1-q_{0,B}).
\end{equation*}

\begin{figure}[th]
\begin{center}
\includegraphics[width=80mm]{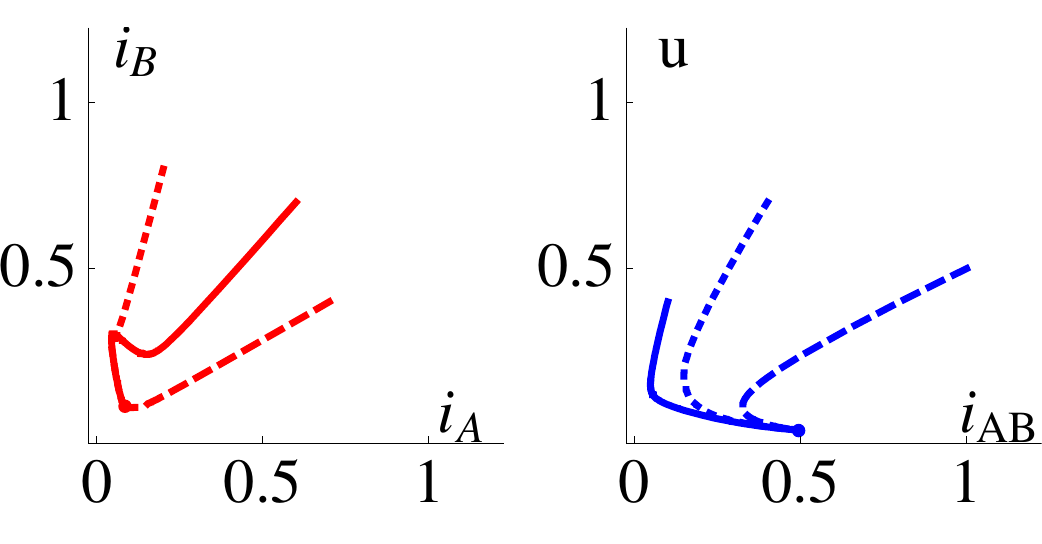}
\caption{With the parameters chosen as described in example
\ref{example_double}, we observe stable coexistence of all three infection
types. Solid, dashed and dotted lines, respectively, show parts of the solution
starting at the same initial values.}\label{fig_6}
\end{center}
\end{figure}

\begin{example}\label{example_double} {\rm Under the above assumptions, let
$q_{A,B}=0.9=q_{B,A}$, the transmission efficacy $\tau_A=\tau_B=0.9$ and the
cost of the infection $\eta_A=\eta_B=1.1$. Then there exists a coexistence
equilibrium of doubly infected and both types of singly infected individuals,
where however the proportion of  doubly infected individuals is much larger
(figure \ref{fig_6}).}
\end{example}

\section{Introduction of age-structure}\label{section:age_structure}

In the previous sections we have seen that infection with \textit{Wolbachia}
gives rise to interesting dynamic behaviour already in unstructured populations.
Clearly, individuals of different ages are subject to different fertility and
mortality rates. We therefore expand our model \eqref{eq_I}--\eqref{eq_U} to
include age-dependent fertility and mortality rates for infected and uninfected
individuals. This leads to nonlinear partial differential equations with
nonlocal boundary conditions that represent the birth
process \citep{Webb,Farkas2006,Farkas2007, Gurtin, Farkas2008}. Although this
results in more complex models, they are still amenable to analytical study.
Here we focus on qualitative questions, when analytical progress is possible; in
particular, how do the stability results for equilibrium solutions compare to
the unstructured case.

Let $i(a,t)$ and $u(a,t)$ denote the densities of infected and uninfected
individuals of age $a$ at time $t$, respectively, where $a\in[0,m]$
(this is not to be confused with the notation in section
\ref{section:single_sex_model}, where they denoted scaled numbers of infected
and uninfected individuals). Then the evolution of the population is governed by
\begin{align}
i_t(a,t)+i_a(a,t)&=-\eta_1(a)(I(t)+U(t))i(a,t), \label{infected_age} \\
u_t(a,t)+u_a(a,t)&=-\eta_2(a)(I(t)+U(t))u(a,t), \label{uninfected_age} \\
i(0,t) &= \tau \int_0^m\beta_1(a)i(a,t)\,\ud a, \label{infected_bdc} \\
u(0,t) &= (1-\tau) \int_0^m\beta_1(a)i(a,t)\,\ud a   \nonumber \\
       &\quad +  \left(1-q\frac{I(t)}{I(t)+U(t)}\right) \int_0^m\beta_2(a)u(a,t)\,\ud a, \label{uninfected_bdc}
\end{align}
where
\begin{equation*}
I(t) = \int_0^m\ i(a,t)\,\ud a, \quad  U(t) = \int_0^m\ u(a,t)\,\ud a,
\end{equation*}
and $\eta_1,\,\eta_2,\,\beta_1$ and $\beta_2$ denote the age-specific mortality
and fertility rates for infected  and uninfected individuals, respectively.
System \eqref{infected_age}--\eqref{uninfected_bdc} is equipped with initial
conditions
\begin{equation*}
i(a,0)=i_0(a),\quad u(a,0)=u_0(a).
\end{equation*}
The parameters $\tau$ and $q$ have the same meaning as in section \ref{section:single_sex_model}.

\subsection{Existence of equilibrium solutions}\label{subsection:existence}
We find the  time-independent solutions of equations \eqref{infected_age} and \eqref{uninfected_age} as
\begin{align}
i_*(a) &= i_*(0)\exp\left\{ -(I_*+U_*)\int_0^a\eta_1(r)\,\ud r \right\}, \label{stateq1} \\
u_*(a) &= u_*(0)\exp\left\{ -(I_*+U_*)\int_0^a\eta_2(r)\,\ud r \right\}, \label{stateq2}
\end{align}
where $i_*(0)$ and $u_*(0)$ satisfy
\begin{align}
i_*(0) &= \tau i_*(0) \int_0^m\beta_1(a) \exp\left\{
-(I_*+U_*)\int_0^a\eta_1(r)\,\ud r \right\}\,\ud a,  \label{stateq3} 
\\
u_*(0) &= (1-\tau) i_*(0) \int_0^m\beta_1(a) \exp\left\{ -(I_*+U_*)\int_0^a\eta_1(r)\,\ud r \right\}\,\ud a \nonumber \\
&\quad + \left(1-q\frac{I_*}{I_*+U_*}\right) u_*(0) \int_0^m\beta_2(a) \exp\left\{ -(I_*+U_*)\int_0^a\eta_2(r)\,\ud r \right\}\,\ud a.\label{stateq4}
\end{align}
Here
\begin{align*}
I_*=\int_0^m i_*(a)\,\ud a, \quad  U_* = \int_0^m u_*(a)\,\ud a.
\end{align*}
First we note that the trivial steady state $(0,0)$ always exists. Next we note
that if $i_*(\cdot)\equiv 0$ then equation \eqref{stateq4} reduces to
\begin{equation}\label{stateq5}
1=\int_0^m\beta_2(a) \exp\left\{ -U_*\int_0^a\eta_2(r)\,\ud r \right\}\,\ud a.
\end{equation}
It then immediately follows from the monotonicity and continuity of the right
hand side of \eqref{stateq5} (as a function of $U_*$) and  the Intermediate
Value Theorem, that a unique disease free equilibrium, given by
\begin{equation*}
u_*(a)=\frac{U_*\,\exp\left\{-U_*\int_0^a\eta_2(r)\,\ud r \right\}}{\int_0^m \exp\left\{-U_*\int_0^a\eta_2(r)\,\ud r \right\}  \ud a},
\end{equation*}
exists if and only
\begin{equation}\label{disfreecond}
\int_0^m\beta_2(a)\ud a>1
\end{equation}
holds. If we look for strictly positive equilibrium solutions $(i_*(a),u_*(a))$ we find that $I_*$ and $U_*$ have to satisfy
\begin{align}
1&=\tau\int_0^m\beta_1(a) \exp\left\{ -(I_*+U_*)\int_0^a\eta_1(r)\,\ud r \right\}\,\ud a,\label{stateq6}  \\
&\!\!\!\! \frac{U_*\left(1-\left(1-q\frac{I_*}{I_*+U_*}\right)\int_0^m\beta_2(a)\exp\left\{-(I_*+U_*)\int_0^a\eta_2(r)\,\ud r\right\}\ud a\right)}
{\int_0^m\exp\left\{-(I_*+U_*)\int_0^a\eta_2(r)\,\ud r\right\}\ud a} \nonumber \\
&=\frac{I_*\left(\tau^{-1}-1\right)}{\int_0^m\exp\left\{-(I_*+U_*)\int_0^a\eta_1(r)\,\ud r\right\}\ud a}.\label{stateq7}
\end{align}
Conversely, if $I_*$ and $U_*$ satisfy equations
\eqref{stateq6}--\eqref{stateq7} then equations
\eqref{stateq1}--\eqref{stateq2} determine uniquely a positive equilibrium
solution. We also see from equation \eqref{stateq6} that
\begin{equation}
\int_0^m\beta_1(a)\ud a>\frac{1}{\tau},\label{stateq8}
\end{equation}
is a necessary condition for the existence of a positive equilibrium. In fact, if equation \eqref{stateq8} holds true, then we can solve equation
\eqref{stateq6} to obtain a unique positive value
\begin{equation}\label{Ustar}
c_1=I_*+U_*.
\end{equation}
A straightforward calculation then leads from equation \eqref{stateq7} to the
following quadratic equation for $I_*$
\begin{equation}\label{Istar}
I_*^2c_3+I_*(1-c_2-c_1c_3+c_4)+c_1c_2-c_1=0, 
\end{equation}
where we have defined 
\begin{align*}
c_2 & = \int_0^m\beta_2(a)\exp\left\{-c_1\int_0^a\eta_2(r)\,\ud r\right\}\ud a,
\quad c_3  = \frac{qc_2}{c_1}, \quad \textrm{and} \\
c_4 & = \frac{(\tau^{-1}-1)\int_0^m\exp\left\{-c_1\int_0^a\eta_2(r)\ud
r\right\}\ud a}{\int_0^m\exp\left\{-c_1\int_0^a\eta_1(r)\ud
r\right\}\ud a}.
\end{align*}
Similarly to the unstructured case, see equation \eqref{quadratic}, we arrive
at a quadratic equation (unless $q=0$) for the infected population size $I_*$.
Of course the calculations now are much more involved since we have
age-dependent fertility and mortality rates. However, for fixed model
ingredients the equilibrium solutions can be determined explicitly, via
equations \eqref{stateq1}--\eqref{stateq2}. In contrast to the unstructured
case, we have necessary conditions on the birth rates for the existence of
non-trivial equilibria.

We summarize our findings in the following theorem.
\begin{theorem}\label{theorem:equilibria}
The equilibrium solutions to equation system
\eqref{infected_age}--\eqref{uninfected_bdc} are given by functions 
\eqref{stateq1}--\eqref{stateq2}  with initial values
\eqref{stateq3}--\eqref{stateq4}, provided that the total populations of
infected and uninfected individuals $I_*$ and $U_*$ given by equations
\eqref{Istar} and \eqref{Ustar} are non-negative. 
\end{theorem}
We note the formal similarity of equations \eqref{Ustar} and \eqref{Istar}
to
the conditions \eqref{u_from_i} and \eqref{quadratic} for the unstructured model
in section \ref{section:single_sex_model}.

\subsection{(In)stability }\label{subsection:stability}
In the previous section we established necessary and sufficient conditions for
the existence of a non-trivial steady-state of system
\eqref{infected_age}--\eqref{uninfected_bdc}. In this section we study stability
properties of the steady states. To this end, first we formally linearise
system \eqref{infected_age}--\eqref{uninfected_bdc} around a steady-state 
solution $(i_*(a),u_*(a))$. We introduce the perturbations
$p(a,t)=i(a,t)-i_*(a)$ and
$s(a,t)=u(a,t)-u_*(a)$ and we use Taylor series expansions of
the fertility and mortality functions. Then we drop the
nonlinear terms to arrive at the linearised system
\begin{align}
p_t(a,t)&+p_a(a,t)=-\eta_1(a)(p(a,t)(I_*+U_*)+i_*(a)(P(t)+S(t))), \label{lin1} \\
s_t(a,t)&+s_a(a,t)=-\eta_2(a)(s(a,t)(I_*+U_*)+u_*(a)(P(t)+S(t))), \label{lin3} \\
p(0,t)&=\tau\int_0^m\beta_1(a)p(a,t)\,\ud a,\label{lin2} \\
s(0,t)&=(1-\tau)\int_0^m\beta_1(a)p(a,t)\,\ud a+\left(1-q\frac{I_*}{I_*+U_*}\right)\int_0^m\beta_2(a)s(a,t)\,\ud a\nonumber \\
&\quad -q\left(\frac{U_*}{(I_*+U_*)^2}P(t)-\frac{I_*}{(I_*+U_*)^2}S(t)\right)\int_0^m\beta_2(a)u_*(a)\,\ud a,\label{lin4}
\end{align}
where
\begin{equation*}
P(t)=\int_0^mp(a,t)\,\ud a,\quad S(t)=\int_0^ms(a,t)\,\ud a.
\end{equation*}
For more detailed calculations we refer the reader to
\citep{Farkas2006,Farkas2007,Farkas2008}, 
where similar type of age and size-structured models where treated. 
It can be shown that the linearised system is governed by a strongly
continuous semigroup of linear operators, which is eventually compact (see
e.g.~\citep{Farkas2007,Farkas2008}). However, this governing semigroup cannot be
shown to be positive, since mortality of both infected and uninfected 
individuals is an increasing function of the total population size. 
Eventual compactness of the governing linear semigroup implies that to study
stability of steady-states it is sufficient to study the point spectrum of the
linearised operator see, e.g.~\citep{Nagel}. The standard way how this can be
carried out is to solve the eigenvalue equation and deduce a characteristic
equation (if possible) whose roots are the eigenvalues of the linearised
operator. We note that the lack of positivity implies that we cannot expect to
establish even local stability results unless the characteristic equation can be
cast in a simple form. We substitute
\begin{equation*}
\begin{pmatrix}
p(a,t) \\ s(a,t)
\end{pmatrix}
=\exp\left\{\lambda t\right\}
\begin{pmatrix}
v(a) \\ w(a)
\end{pmatrix}
\end{equation*}
into the linearised equations \eqref{lin1}--\eqref{lin4}. This yields
\begin{align}
v'(a) & = -v(a)\left(\lambda+\eta_1(a)(I_*+U_*)\right)-\eta_1(a)i_*(a)(V+W),
\label{lin5} \\
w'(a) & =
-w(a)\left(\lambda+\eta_2(a)(I_*+U_*)\right)-\eta_2(a)u_*(a)(V+W), \label{lin7}
\\
v(0) & = \tau\int_0^m\beta_1(a)v(a)\,\ud a, \label{lin6} \\
w(0) & = (1-\tau)\int_0^m\beta_1(a)v(a)\,\ud
a+\left(1-q\frac{I_*}{I_*+U_*}\right)\int_0^m\beta_2(a)w(a)\,\ud a \nonumber \\
& \quad +q\frac{I_*W-U_*V}{(I_*+U_*)^2}\int_0^m\beta_2(a)u_*(a)\,\ud a,\label{lin8}
\end{align}
where
\begin{equation*}
V=\int_0^mv(a)\,\ud a,\quad\quad W=\int_0^mw(a)\,\ud a.
\end{equation*}
Hence $\lambda\in\mathbb{C}$ is an eigenvalue if and only if the nonlocal system
\eqref{lin5}--\eqref{lin8} admits a non-trivial solution. The solution of the
differential equations \eqref{lin5} and \eqref{lin7} is
\begin{align}
v(a) &
=f_{\lambda}^1(a)\left(v(0)-\int_0^a\frac{\eta_1(x)i_*(x)(V+W)}{f_\lambda^1(x)}\,
\ud x\right),\label{lin9} \\
w(a) &
=f_\lambda^2(a)\left(w(0)-\int_0^a\frac{
\eta_2(x)u_*(x)(V+W) } { f_\lambda^2(x) } \, \ud x\right),\label{lin10}
\end{align}
where we have introduced
\begin{equation*}
f_\lambda^i(a)=\exp\left\{-\int_0^a\lambda+\eta_i(y)(I_*+U_*)\,\ud
y\right\},\quad i=1,2.
\end{equation*}
Next we substitute the solutions \eqref{lin9} and \eqref{lin10} into the
boundary conditions \eqref{lin6} and \eqref{lin8}
 and integrate the solution \eqref{lin9} and \eqref{lin10} from zero to $m$ to
arrive at a four-dimensional homogeneous system for 
the unknowns $v(0), w(0), V$ and $W$. This homogeneous system admits a
non-trivial solution 
if and only if the determinant of the coefficient matrix equals zero. We can formulate the following theorem.
\begin{theorem}\label{theorem:eigenvalues}
$\lambda$ is an eigenvalue of the linearised operator if and only if it satisfies the equation
\begin{equation}\label{matrix}
\begin{aligned}
K(\lambda)&= \det
\left( \begin{array}{llll}
\tau a_5(\lambda) - 1 & 0 & -\tau a_6(\lambda) & -\tau a_6(\lambda)\\
(1-\tau)a_5(\lambda) & a_7(\lambda)-\frac{q\,I_*\,a_7(\lambda)}{I_*+U_*}-1  & a_8(\lambda) & a_{10}(\lambda) \\
a_1(\lambda) & 0 & -a_2(\lambda) -1 & -a_2(\lambda)\\
0 & a_3(\lambda) & -a_4(\lambda) & -a_4(\lambda)-1
\end{array}\right) \\ &=0,
\end{aligned}
\end{equation}
where
\begin{align*}
a_1(\lambda)&= \int_0^mf_\lambda^1(a)\,\ud a,\quad\quad\quad\,
a_2(\lambda)=\int_0^mf_\lambda^1(a)\int_0^a
\frac{\eta_1(x)i_*(x)}{f_\lambda^1(x)}\,\ud x\,\ud a, \\
a_3(\lambda)&=  \int_0^mf_\lambda^2(a)\,\ud a,\quad\quad\quad\,
a_4(\lambda)=\int_0^mf_\lambda^2(a)\int_0^a
\frac{\eta_2(x)u_*(x)}{f_\lambda^2(x)}\,\ud x\,\ud a, \\
a_5(\lambda)&=  \int_0^m\beta_1(a)f_\lambda^1(a)\,\ud a,\quad
a_6(\lambda)=\int_0^m\beta_1(a)f_\lambda^1(a)\int_0^a\frac{\eta_1(x)i_*(x)}{
f_\lambda^1(x)}\,\ud x\,\ud a, \\
a_7(\lambda)&=  \int_0^m\beta_2(a)f_\lambda^2(a)\,\ud a,\quad
a_9(\lambda)=\int_0^m\beta_2(a)f_\lambda^2(a)\int_0^a\frac{\eta_2(x)u_*(x)}{
f_\lambda^2(x)}\,\ud x\,\ud a, \\
a_8(\lambda)&=
(\tau-1)a_6(\lambda)+\left(\frac{q\,I_*}{I_*+U_*}-1\right)a_9(\lambda)-\frac{q\,
U_*}{(I_*+U_*)^2}\int_0^m\beta_2(a)u_*(a)\,\ud a, \\
a_{10}(\lambda)&=
(\tau-1)a_6(\lambda)+\left(\frac{q\,I_*}{I_*+U_*}-1\right)a_9(\lambda)+\frac{q\,
I_*}{(I_*+U_*)^2}\int_0^m\beta_2(a)u_*(a)\,\ud a.
\end{align*}
\end{theorem}
\noindent It follows from the growth behaviour of the functions $f^i_\lambda$ that 
\begin{equation}
    \lim_{\lambda\rightarrow +\infty} K(\lambda)= \det
\begin{pmatrix} -1 & 0 & 0 & 0\\
                0 & -1 & C_1 & C_2\\
                0 & 0 & -1 & 0\\
                0 & 0 & 0 & -1\\
\end{pmatrix} = 1,
\end{equation}
the limit being taken in $\mathbb R$, and $C_1,C_2$  are constants. Hence we can
formulate the following general instability criterion, which follows immediately
from the Intermediate Value Theorem.

\begin{theorem}\label{instabeig}
The stationary solution $(i_*(a),u_*(a))$ of
equations \eqref{infected_age}--\eqref{uninfected_bdc} is unstable if $K(0)<0$.
\end{theorem}
As we can see the characteristic equation \eqref{matrix} is rather complicated in ge\-ne\-ral, hence we only consider some interesting special cases 
when analytical progress is possible.

\subsubsection{The trivial steady state}

We consider the stability of the steady state $i_*\equiv 0,\,u_*\equiv 0$. Note that in this case the characteristic equation \eqref{matrix} 
reduces to
\begin{equation}\label{matrix2}
K(\lambda)= \det
\left( \begin{array}{llll}
\tau a_5(\lambda) - 1 & 0 & 0 & 0\\
(1-\tau)a_5(\lambda) & a_7(\lambda)-1  & 0 & 0 \\
a_1(\lambda) & 0 & -1 & 0 \\
0 & a_3(\lambda) & 0 & -1
\end{array}\right)=0,
\end{equation}
which leads to the equation
\begin{equation}\label{trivialchareq}
(\tau a_5(\lambda)-1)(a_7(\lambda)-1)=0.
\end{equation}
Therefore, $\lambda\in \mathbb{C}$ is an eigenvalue if and only if $\lambda$
satisfies either of the two equations
\begin{equation}
1=\tau\int_0^m\beta_1(a)e^{-\lambda a}\,\ud a,\quad 1=\int_0^m\beta_2(a)e^{-\lambda a}\,\ud a.
\end{equation}
We can formulate the following theorem.
\begin{theorem}
The trivial steady state is locally asymptotically stable if 
\begin{equation}
\tau\int_0^m\beta_1(a)\,\ud a <1\quad \text{and}\quad \int_0^m\beta_2(a)\,\ud a<1.
\end{equation}
On the other hand, if either
\begin{equation}
\tau\int_0^m\beta_1(a)\,\ud a >1\quad \text{or} \quad \int_0^m\beta_2(a)\,\ud a>1.
\end{equation}
holds, then the trivial steady state is unstable.
\end{theorem}

\subsubsection{The disease free steady state}

Consider the disease free steady state, i.e. $i_*\equiv 0$, which exists by
condition \eqref{disfreecond} if and only if $\int_0^m\beta_2(a)\,\ud a>1$. 
In this case the characteristic equation \eqref{matrix} can be written as
\begin{equation}\label{matrix3}
\begin{aligned}
K(\lambda)&=  (\tau a_5(\lambda)-1) \\
&\phantom{=}\times\det
\left( \begin{array}{lll}
 a_7(\lambda)-1 & -a_9(\lambda)-\frac{q\int_0^m\beta_2(a)u_*(a)\,\ud a}{U_*} & -a_9(\lambda)\\
0  & -1 & 0 \\
a_3(\lambda) & -a_4(\lambda) & -a_4(\lambda)-1 
\end{array}\right) \\
&=  (\tau a_5(\lambda)-1)\left[
(-a_4(\lambda)-1)(a_7(\lambda)-1)+a_3(\lambda)a_9(\lambda)\right] \\
&=  0.
\end{aligned}
\end{equation}
This again splits into two equations. The first one is easy to analyse, since it
can be written as
\begin{equation}
1=\tau\int_0^m\beta_1(a)\exp\left\{-U_*\int_0^a\eta_1(x)\,\ud x\right\}e^{-\lambda a}\,\ud a.
\end{equation}
\begin{theorem}\label{disfreeunstable}
If 
\begin{equation}
\tau\int_0^m\beta_1(a)\exp\left\{-U_*\int_0^a\eta_1(x)\,\ud x\right\}\,\ud a>1,\label{instabconddisfree}
\end{equation}
where $U^*$ satisfies \eqref{stateq5}, then the disease free steady state is
unstable.
\end{theorem}
\begin{remark}
Provided that equation  \eqref{matrix3} has a dominant real solution 
$\lambda$, it is shown that condition \eqref{instabconddisfree} in Theorem \ref{disfreeunstable} is indeed necessary and sufficient 
for the instability of the disease free steady state. However, as we have noted before the governing linear semigroup cannot shown to be positive, 
hence we cannot establish the existence of a dominant real root of the
characteristic function \eqref{matrix}.
\end{remark}

\subsubsection{Complete transmission of disease}

In case of complete transmission of disease, i.e. when $\tau=1$, equations
\eqref{infected_age}--\eqref{uninfected_bdc} can be written in 
the following form
\begin{align}
i_t(a,t)+i_a(a,t)&=-\mu_1(a,I(t),U(t))i(a,t), \label{infected_age2} \\
u_t(a,t)+u_a(a,t)&=-\mu_2(a,I(t),U(t))u(a,t), \label{uninfected_age2} \\
i(0,t) &= \int_0^m\beta_1(a)i(a,t)\,\ud a, \label{infected_bdc2} \\
u(0,t) &= \int_0^m\beta^q_2(a,I(t),U(t))u(a,t)\,\ud a, \label{uninfected_bdc2}
\end{align}
where 
\begin{align}
 \mu_i(a,I(t),U(t))&=\eta_i(a)(I(t)+U(t)),\quad i=1,2,
\label{mortality} \\
\beta_2^q(a,I(t),U(t))&=\beta_2(a)\left(1-q\frac{I(t)}{I(t)+U(t)}\right).
\label{fertility}
\end{align}
Hence model \eqref{infected_age2}--\eqref{uninfected_bdc2} is a special case of
the $n$-species age-structured system considered in \citep{Farkas2006}, where
the coupling occurs due to competition for resources and due to the inhibition
of the proliferation of the uninfected population.
 
In \citep{Farkas2006} we deduced a very general instability condition, which we
recall for the case $n=2$ for the readers convenience (see Theorem 2.3 in
\citep{Farkas2006}).
\begin{theorem}\label{2species}
A strictly positive stationary solution
$(i_*,u_*)\in (\RR_{>0})^2$ of \eqref{infected_age2}--\eqref{uninfected_bdc2} is
unstable if the partial derivatives of the net reproduction rates of the
infected, respectively uninfected populations satisfy
\begin{equation}\label{instabilitycond}
R^1_{I}(I_*,U_*)R^2_{U}(I_*,U_*)-R^1_{U}(I_*,U_*)R^2_{I}(I_*,U_*)<0,
\end{equation}
\end{theorem}
Taking into account \eqref{mortality}--\eqref{fertility} we have
\begin{align}
R^1(I,U)= & \int_0^m\beta_1(a)\exp\left\{-(I+U)\int_0^a\eta_1(x)\,\ud x\right\}\,\ud a,\label{rep1} \\
R^2(I,U)= & \int_0^m\beta_2(a)\left(1-q\frac{I}{I+U}\right)\exp\left\{-(I+U)\int_0^a\eta_2(x)\,\ud x\right\}\,\ud a\label{rep2}.
\end{align}
From equation \eqref{rep1} we obtain
\begin{align*}
&R_{I}^1(I_*,U_*)=R_{U}^1(I_*,U_*)\\
& =-\int_0^m\beta_1(a)\left(\int_0^a\eta_1(x)\,\ud
x\right)\exp\left\{-(I_*+U_*)\int_0^a\eta_1(x)\,\ud x\right\}\,\ud
a<0,
\end{align*}
unless $\beta_1\equiv 0$ or $\eta_1\equiv 0$. Also from equation \eqref{rep2} we obtain
\begin{align*}
& R^2_{I}(I_*,U_*)= \\
&
-\int_0^m\beta_2(a)\left(1-q\frac{I_*}{I_*+U_*}\right)\left(\int_0^m\eta_2(x)\,
\ud x\right)\exp\left\{-(I_*+U_*)\int_0^a\eta_2(x)\,\ud x\right\}\,\ud a \\
& -\int_0^m\beta_2(a)\exp\left\{-(I_*+U_*)\int_0^a\eta_2(x)\,\ud
x\right\}q\frac{U_*}{(I_*+U_*)^2}\,\ud a,\\
& R^2_{U}(I_*,U_*)= \\
&
-\int_0^m\beta_2(a)\left(1-q\frac{I_*}{I_*+U_*}\right)\left(\int_0^m\eta_2(x)\,
\ud x\right)\exp\left\{-(I_*+U_*)\int_0^a\eta_2(x)\,\ud x\right\}\,\ud a \\
& +\int_0^m\beta_2(a)\exp\left\{-(I_*+U_*)\int_0^a\eta_2(x)\,\ud
x\right\}q\frac{I_*}{(I_*+U_*)^2}\,\ud a.
\end{align*}
Hence $R^2_{U}(I_*,U_*)>R^2_{I}(I_*,U_*)$ for every 
strictly positive steady state, unless $\beta_2\equiv 0$ or
$q=0$ (in which case $0$ is the strictly dominant eigenvalue of the linearised
operator). We summarize our findings in the following theorem.
\begin{theorem}\label{compltrans}
Assume that  $\tau=1$, $q\ne 0$, and $\beta_1,\beta_2,\eta_1$ are not
identically zero. Then any strictly positive
 steady state of equations \eqref{infected_age}--\eqref{uninfected_bdc} 
is unstable.
\end{theorem}
In other words,  there is no coexistence of infected and
uninfected populations, This corresponds to the instability of the equilibrium
solution $(i_1,u_1)$ in the left panel of figure \ref{fig_2} for the
unstructured case.

\section{Discussion}\label{section:discussion}
In the present work we introduced and studied differential equation models for
the dynamics of populations infected with \textit{Wolbachia}. First we built
ordinary differential equation models, in which we have implemented fitness
costs of an infection as increased mortalities while keeping the birth rates
equal for all infection statuses. It is equally appropriate to reduce birth
rates for infected individuals and (for the sake of simplicity) then to assign
the same mortality to all individuals. This leads to the following model for the
case of a single \textit{Wolbachia} in an asexual population
\begin{align*}
\frac{di}{dt} &= (\mu\tau -(i+u))i, \\
\frac{du}{dt} &= \mu(1-\tau)i +\left( 1-q\frac{i}{i+u}  -(i+u)\right)u,
\end{align*}
where $\mu\in[0,1]$ is the reduced fecundity of infected individuals. This
results in similar formulas for equilibrium solutions as \eqref{coexist} and
vector fields as in figure \ref{fig_2}. An experimentally testable prediction of
our model \eqref{eq_I}--\eqref{eq_U} is that there are no persistent
\textit{Wolbachia} strains with a transmission efficacy less than $\frac{3}{4}$
(see the region for existence of the observed stable equilibrium
$(u_2,i_2)$ in figure  \ref{fig_1}).

Our model for multiple infections is novel insofar it allows the theoretical
biologist to adapt it to a case of special interest (with or without doubly
infected individuals, with or without mutual incompatibility). This should help
to gain a more unified perspective than what was possible from models created
for each purpose individually. In the case of mutual compatibility we saw that
strains with higher transmission efficacy or lower mortality due to infection
establish themselves over competitors. This is in good agreement with other
predictions from discrete population genetics models
\citep{Turelli1994,Haygood2009}. Although the model for infections with multiple
mutually incompatible  strains in section \ref{section:multiple_strains} is too
complicated for all of its equilibrium solutions to be written down explicitly,
it can be analyzed to a certain degree by identifying invariant subspaces. By
numerical simulations we provided evidence for the absence of coexistence of
singly infected individuals, apart for exceptional choices of parameters and
initial values. This situation changes, if doubly infected individuals are 
present that can lose  one of their strains when giving birth.

We have expanded the simple ordinary differential equation model from section
\ref{section:single_sex_model} by introducing age-structure. Of the extensive
literature about structured populations, let us mention the monographs
\citep{Cushing,Webb, Metz}, the classical paper \citep{Gurtin} and the recent
collection \citep{MagalRuan}. Clearly, age-structured models allow a much finer
level of detail to be incorporated, but also pose greater analytical challenges.
 Nevertheless, we have shown that existence of equilibrium solutions  and their
stability properties can be investigated in a straightforward fashion. We saw
that unlike in the unstructured case, existence of a disease free equilibrium is
now subject to a condition on the integral of the birth rate. We also obtained
analytical results in some special cases which allow at least a partial
characterization of the dynamic behaviour of the system once the model
parameters are fixed.

\textit{Wolbachia} together  with \textit{Cardinium} are the two bacterial
infections of arthropods that cause cytoplasmic incompatibility. In this work we
have focused on the case of diplodiploid species where cytoplasmic
incompatibility results, with a certain probability, in embryonic death. In
future work, we will include separate sexes into our models, see the book by
Iannelli \textit{et al.}~\citep{Iannelli} for a comprehensive introduction to
gender-structured populations. Then it will be possible to study gender-specific
effects of the \textit{Wolbachia} infection. These become even more important in
haplodiploid organisms such as bees, ants and wasps, where cytoplasmic
incompatibility is vastly more complex \citep{Vautrin2007,Stouthamer1997}
(male-development, thelytokous parthenogenesis etc.).

\section*{\normalsize{Acknowledgments.}}
J.~Z. Farkas is thankful to the Centre de Recerca Mathem\`{a}tica, Universitat
Aut\`{o}noma de Barcelona for their hospitality while being a participant in the
research programme "Mathematical Biology: Modelling and Differential Equations"
during 01/2009-06/2009. J.~Z. Farkas was also supported by a personal research
grant from the Carnegie Trust for the Universities of Scotland. P.~Hinow was
supported partly by the NSF through an IMA postdoctoral fellowship. Part
of the work on this paper was done while P.~Hinow visited the Department of
Computing Science and Mathematics at the University of Stirling. Financial
support from the Edinburgh Mathematical Society during this visit is greatly
appreciated. We thank Jan Engelst\"adter (Institute for Integrative Biology,
Swiss Federal Institute of Technology, Z\"urich, Switzerland) for helpful
remarks and advice on literature. We are much indebted to the anonymous referees
for their helpful comments.
\bibliography{wolbachia}
\end{document}